\documentclass[]{amsart}
\usepackage{dcolumn}
\usepackage{amsmath}

\usepackage{graphicx}
\usepackage{fancyvrb}
\usepackage{bbold}
\newtheorem{defn}{Definition}
\newtheorem{conj}{Conjecture}

\usepackage{color}
\definecolor{brown}{rgb}{0.8,0.6,0.3}
\definecolor{dgreen}{rgb}{0.2,0.4,0.3}

\usepackage[pdftex,pdfauthor={Richard J. Mathar}, colorlinks=true,urlcolor=blue,linkcolor=red,citecolor=brown,
pdfkeywords={{Combinatorics},{Tiling},{Squares}},bookmarksopen=true,pdfpagelayout=OneColumn]{hyperref}

\begin{document}

\title[Tiling rectangles with monomers and squares]{Tiling $n\times m$ rectangles with $1\times 1$ and $s\times s$ squares}

\author{Richard J. Mathar}
\urladdr{http://www.mpia.de/~mathar}
\address{K\"onigstuhl 17, 69117 Heidelberg, Germany}

\date{\today}

\begin{abstract}
We consider tilings 
of a rectangle which is $n$ units wide and $m$ units long by
non-overlapping $1\times 1$ squares and $s\times s$ squares.
Bivariate generating functions are computed with the Transfer Matrix Method
for moderately large
but fixed widths $n$ as a function of the parameter $m$ and of the number
of $s\times s$ squares in the rectangle.
\end{abstract}

\maketitle

\section{Definitions}
We consider the combinatorial problem of placing non-overlapping
squares of shape $s\times s$
into rectangles of shape $n\times m$. Comparing the areas 
$nm$ of the hosting rectangle and the area $s^2$
of the individual square we find a trivial upper limit for the number $k$
of $s\times s$ squares that fit into the rectangle:
\begin{equation}
0\le k \le nm/s^2 .
\end{equation}
The area in the rectangle that is not covered
by the $s\times s$ squares is tiled
with
$1\times 1$ squares (monomers), of which there are $nm-ks^2$.

The manuscript is basically an industrial scale evaluation of
Heubach's tilings \cite{HeubachCN140}.

\begin{figure}
\includegraphics[scale=0.7]{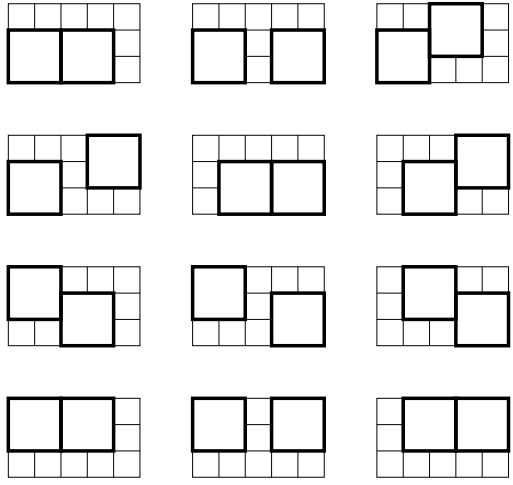}
\caption{Illustration of the $T_{3\times 5}(2,2)=12$ ways of placing
$k=2$ squares of edge length $s=2$ in a $n\times m = 3\times 5$ rectangle.
}
\label{fig.exa}
\end{figure}

\begin{defn}
$T_{n\times m}(s,k)$ is the number of ways of tiling the $n\times m$ rectangle
with $k$ non-overlapping squares of shape $s\times s$ and 
with $nm-ks^2$ unit squares.
 Distributions obtained by flipping or rotating
the rectangle are considered distinct and counted with multiplicity.
\end{defn}
A basic example of such counting on commensurate grids are the tilings of Figure \ref{fig.exa},
which shows all variants of distributing two $2\times 2$ squares on a $3\times 5$ board.
The total number of geometries that does not resolve how many
squares fill the rectangle is:
\begin{defn}
The number of ways of tiling the $n\times m$ rectangle
with $1\times 1$ and $s\times s$ squares is
\begin{equation}
T_{n\times m}(s)\equiv \sum_{k=0}^{\lfloor nm/s^2\rfloor} 
T_{n\times m}(s,k).
\label{eq.rsum}
\end{equation}
\end{defn}
 
Bivariate ordinary generating functions will be noted as follows:
\begin{defn}
\begin{equation}
\sum_{m\ge 0} \sum_{k\ge 0} T_{n\times m}(s,k) z^m t^k \equiv T_n(s,z,t).
\label{eq.2g}
\end{equation}
\end{defn}

\section{Symmetries}
There is always one way of filling the rectangle with monomers only:
\begin{equation}
T_{n\times m}(s,0)=1.
\end{equation}
If only one $s\times s$ square is to be placed, it can be rooted at any
of the $n-s+1$ vertical and $m-s+1$ horizontal grid points:
\begin{equation}
T_{n\times m}(s,1)=(n-s+1)(m-s+1),\quad m,n\ge s. 
\label{eq.k1}
\end{equation}

Rotating the rectangle (and rotating the embedded squares with it) does
not change the count; it is symmetric with respect to interchange
of $n$ and $m$:
\begin{equation}
T_{n\times m}(s,k)=T_{m\times n}(s,k).
\end{equation}

If the rectangle has a width and height that are integer multiples $is$
and $js$
of the square edge $s$, there is one configuration with full coverage
by the $s\times s$ squares, not using interstitial monomers:
\begin{equation}
T_{is\times js}(s,ij) = 1,\quad i,j\ge 1.
\label{eq.ij}
\end{equation}

If the width or height are too small, only the $1\times 1$ squares fit in:
\begin{equation}
T_{n\times m}(s,k) = \delta_{0,k}, \quad (s>m \vee s> n).
\end{equation}

If the width equals the square size, $s=n$, there is an obvious bijection
to filling a line
with monomers and straight $s$-ominos, equal to the number of compositions
of $m$ into 1's and $k$ $s$'s:
\begin{equation}
T_{s\times m}(s,k) = \binom{m-(s-1)k}{k}.
\label{eq.sonce}
\end{equation}

\section{Transfer Matrix Technique}
The bivariate generating functions are constructed with
a variant of the Transfer Matrix technique specialized
to the tiling combinatorics \cite{MatharArxiv1406}.
The growth front of incrementally adding $s\times s$ squares
or leaving free space (that is, adding  $1\times 1$ squares)
is encoded in an integer vector of $n$ entries, one per ``lane.''
These vectors are states in the directed graph of all possible fronts,
and we construct all states that are reachable starting from a vector
of all-zeros by repeatedly/recursively trying to attach a set of
squares to the front to reach the next state.
The only difference to the earlier strategy
with univariate generating functions
\cite{MatharArxiv1406} is that a transition between two states
does not only introduce a factor $z$ to indicate that the
base line is rolled up by one unit, but also another factor $t^k$
where $k$ is the number of blocks of width $s$ that recede
to the back of the front line.

The implementation of this is put into concrete in the appendix.
A Java program constructs all the reachable states, counts them
to get the size of the Transfer Matrix, fills the matrix with
the factors of either zero (if the row state is not reachable directly
from the column state) or $zt^k$, and writes a Maple program that
actually solves the linear system of equations to get the head
element of the inverse. The mere limitation is the
patience needed to execute the Maple program if the Transfer Matrix
(node count in the digraph) exceeds a dimension of, say, $350$.

\section{Results}
The results are tabulated in the following format of a 4-dimensional
table:
Each line contains $s$, then $n$, then $m$, then a colon,
then a sequence of $T_{n\times m}(s,k)$ for $k=0,1,2,\ldots$
and finally another colon and the row sum (\ref{eq.rsum}).
Univariate generating functions of the row sums are obtained by inserting 
$t=1$ into the bivariate generating function (\ref{eq.2g}).

\subsection{$2\times 2$ Squares}
Putting $2\times 2$ squares into rectangles yields:
\small
\begin{Verbatim}
 2 1 1:  1 : 1
 2 2 1:  1 : 1
 2 2 2:  1 1 : 2
 2 3 1:  1 : 1
 2 3 2:  1 2 : 3
 2 3 3:  1 4 0 : 5
 2 4 1:  1 0 : 1
 2 4 2:  1 3 1 : 5
 2 4 3:  1 6 4 0 : 11
 2 4 4:  1 9 16 8 1 : 35
 2 5 1:  1 0 : 1
 2 5 2:  1 4 3 : 8
 2 5 3:  1 8 12 0 : 21
 2 5 4:  1 12 37 34 9 0 : 93
 2 5 5:  1 16 78 140 79 0 0 : 314
\end{Verbatim}
\normalsize
The row sums are tabulated in \cite[A245013]{EIS}.
If the rectangle is only 2 units wide, the problem is equivalent
to counting monomer-dimer coverings of a stripe, and the
Fibonacci numbers appear as row sums \cite[A011973]{EIS}:

\small
\begin{Verbatim}
 2 2 1:  1 : 1
 2 2 2:  1 1 : 2
 2 2 3:  1 2 : 3
 2 2 4:  1 3 1 : 5
 2 2 5:  1 4 3 : 8
 2 2 6:  1 5 6 1 : 13
 2 2 7:  1 6 10 4 : 21
\end{Verbatim}
\normalsize
The generating function is
\begin{equation}
T_2(2,z,t) = \frac{1}{1-z-z^2t}.
\end{equation}

If the rectangle is 3 units wide, each of the $2\times 2$
squares has one more place to slide sideways, 
$T_{3\times m}(2,k)= 2^k T_{2\times m}(2,k)$ \cite[A128099]{EIS}:
\small
\begin{Verbatim}
 2 3 1:  1 : 1
 2 3 2:  1 2 : 3
 2 3 3:  1 4 0 : 5
 2 3 4:  1 6 4 0 : 11
 2 3 5:  1 8 12 0 : 21
 2 3 6:  1 10 24 8 0 : 43
 2 3 7:  1 12 40 32 0 0 : 85
 2 3 8:  1 14 60 80 16 0 0 : 171
 2 3 9:  1 16 84 160 80 0 0 : 341
 2 3 10:  1 18 112 280 240 32 0 0 : 683
\end{Verbatim}
\normalsize
The row sums are \cite[A001045]{EIS}
\begin{equation}
T_{3\times m}(2) =\frac{2^{n+1}+(-)^n}{3}.
\end{equation}
The generating function is
\begin{equation}
T_3(2,z,t) = \frac{1}{1-z-2z^2t}.
\end{equation}

If the rectangle is 4 units wide, there are options 
to stack the $2\times 2$ squares \cite[A128101]{EIS}:
\small
\begin{Verbatim}
 2 4 1:  1 0 : 1
 2 4 2:  1 3 1 : 5
 2 4 3:  1 6 4 0 : 11
 2 4 4:  1 9 16 8 1 : 35
 2 4 5:  1 12 37 34 9 0 : 93
 2 4 6:  1 15 67 105 65 15 1 : 269
 2 4 7:  1 18 106 248 250 108 16 0 : 747
 2 4 8:  1 21 154 490 726 522 176 24 1 : 2115
 2 4 9:  1 24 211 858 1736 1824 994 260 25 0 : 5933
 2 4 10:  1 27 277 1379 3604 5148 4090 1770 385 35 1 : 16717
 2 4 11:  1 30 352 2080 6735 12438 13406 8424 2971 530 36 0 : 47003
 2 4 12:  1 33 436 2988 11615 26691 37150 31598 16207 4787 736 48 1 : 132291
\end{Verbatim}
\normalsize
The generating function is
\begin{verbatim}
T_4(2,z,t) = (-z*t+1) / (-z*t -2*z^2*t +z^3*t^2 +z^3*t^3 +1 -z -z^2*t^2),
\end{verbatim}
and for the row sums \cite[A054854]{EIS}
\begin{verbatim}
T_4(2,z,1) = (1-z) / (-2*z -3*z^2 +2*z^3 +1).
\end{verbatim}

If the rectangle is 5 units wide \cite[A054855]{EIS},
\small
\begin{Verbatim}
 2 5 1:  1 0 : 1
 2 5 2:  1 4 3 : 8
 2 5 3:  1 8 12 0 : 21
 2 5 4:  1 12 37 34 9 0 : 93
 2 5 5:  1 16 78 140 79 0 0 : 314
 2 5 6:  1 20 135 382 454 194 27 0 : 1213
 2 5 7:  1 24 208 824 1566 1344 408 0 0 : 4375
 2 5 8:  1 28 297 1530 4103 5670 3698 926 81 0 0 : 16334
 2 5 9:  1 32 402 2564 9009 17696 18738 9636 1847 0 0 0 : 59925
 2 5 10:  1 36 523 3990 17484 45274 68545 57648 24067 3988 243 0 0 : 221799
 2 5 11:  1 40 660 5872 30984 100608 201832 244080 167565 57940 7698 0 0 0 : 817280
 2 5 12:  1 44 813 8274 51221 201278 508624 818722 812110 464710 135715 16060 729 0 0 0 : 3018301
\end{Verbatim}
\normalsize
with generating function
\begin{verbatim}
T_5(2,z,t) = (-z^2*t^2 -z*t +1) / (-3*z^2*t -z -t^2*z^3 +3*z^4*t^4 
                                  +3*z^3*t^3 -4*z^2*t^2 -z*t +1)
\end{verbatim}
and with row sums
\begin{verbatim}
T_5(2,z,1) = (-z^2-z+1) / (-7*z^2 -2*z +2*z^3 +3*z^4 +1).
\end{verbatim}

If the rectangle is 6 units wide \cite[A063650]{EIS}
\small
\begin{Verbatim}
 2 6 1:  1 0 : 1
 2 6 2:  1 5 6 1 : 13
 2 6 3:  1 10 24 8 0 : 43
 2 6 4:  1 15 67 105 65 15 1 : 269
 2 6 5:  1 20 135 382 454 194 27 0 : 1213
 2 6 6:  1 25 228 964 1987 1974 978 242 27 1 : 6427
 2 6 7:  1 30 346 1976 6014 9856 8544 3760 796 64 0 : 31387
 2 6 8:  1 35 489 3543 14510 34475 47394 37282 16882 4378 619 42 1 : 159651
 2 6 9:  1 40 657 5790 30075 95466 186715 223696 162531 70502 17673 2340 125 0 : 795611
\end{Verbatim}
\normalsize
with generating function
\begin{verbatim}
T_6(2,z,t) = ( -z^5*t^7  +z^4*t^6 +3*z^3*t^4 +2*z^3*t^3 -2*z^2*t^3 -2*z^2*t^2 -2*z*t +1)
 /(z^7*t^10 +3*z^7*t^9 -z^6*t^9 -3*z^6*t^8 +z^6*t^7 -4*z^5*t^7 -15*z^5*t^6 
+3*z^4*t^6 -9*z^5*t^5 +11*z^4*t^5 +12*z^4*t^4 +5*z^3*t^4 +2*z^4*t^3 +10*z^3*t^3 
-3*z^2*t^3 -8*z^2*t^2 -3*z^2*t -2*z*t +1 -z),
\end{verbatim}
and with row sums
\begin{verbatim}
T_6(2,z,1) = (z^4 -5*z^2 -z +1)/( -4*z^6 -z^5 +27*z^4 -z^3 -16*z^2 -2*z +1).
\end{verbatim}
If it is 7 units wide
\small
\begin{Verbatim}
 2 7 1:  1 0 : 1
 2 7 2:  1 6 10 4 : 21
 2 7 3:  1 12 40 32 0 0 : 85
 2 7 4:  1 18 106 248 250 108 16 0 : 747
 2 7 5:  1 24 208 824 1566 1344 408 0 0 : 4375
 2 7 6:  1 30 346 1976 6014 9856 8544 3760 796 64 0 : 31387
 2 7 7:  1 36 520 3920 16834 42368 62266 51504 21792 3600 0 0 0 : 202841
\end{Verbatim}
\normalsize
with generating functions
\begin{verbatim}
T_7(2,z,t) = ( -6*z^6*t^9 +14*z^4*t^6 +1 -2*z*t -5*z^2*t^2 -7*z^2*t^3 +8*z^3*t^4 
+3*z^4*t^5 +z^3*t^3 -6*z^5*t^7) / (24*z^8*t^12 +8*z^8*t^11 +24*z^7*t^10 +18*z^7*t^9 
-62*z^6*t^9 -48*z^6*t^8 -6*z^6*t^7 -38*z^5*t^7 -42*z^5*t^6 +42*z^4*t^6 -5*t^5*z^5 
+65*z^4*t^5 +30*z^4*t^4 +16*z^3*t^4 +t^3*z^4 +4*z^3*t^3 -11*z^2*t^3 -3*t^2*z^3 
-15*z^2*t^2 -4*t*z^2 -2*z*t +1 -z)
\end{verbatim}
and
\begin{verbatim}
T_7(2,z,1) = ( -6*z^6 +17*z^4 +1 -2*z -12*z^2 +9*z^3 -6*z^5) / 
(32*z^8 +42*z^7 -116*z^6 -85*z^5 +138*z^4 +17*z^3 -30*z^2 -3*z +1).
\end{verbatim}
If it is 8 units wide
\small
\begin{Verbatim}
2 8 0: 1 : 1
2 8 1: 1 0 0 : 1
2 8 2: 1 7 15 10 1 : 34
2 8 3: 1 14 60 80 16 0 0 : 171
2 8 4: 1 21 154 490 726 522 176 24 1 : 2115
2 8 5: 1 28 297 1530 4103 5670 3698 926 81 0 0 : 16334
2 8 6: 1 35 489 3543 14510 34475 47394 37282 16882 4378 619 42 1 : 159651
2 8 7: 1 42 730 6872 38562 134088 291908 394202 321408 150950 38296 4944 256 0 0 : 1382259
\end{Verbatim}
\normalsize
with generating functions
\begin{verbatim}
T_8(2,z,t) = (1 -85*z^9*t^16 -53*z^8*t^14 +2*z^7*t^14 +49*t^6*z^4 +12*z^6*t^11
   +206*t^12*z^7 -10*z^9*t^17 +2*z^6*t^12 -45*z^8*t^13 +43*z^7*t^13 -16*z^8*t^15
   +271*z^7*t^11 +11*z^4*t^5 -41*t^9*z^6 -13*z^2*t^3 +5*t^3*z^3 +6*z^3*t^6 -2*z*t^2
   -3*z*t +40*z^3*t^4 +33*z^10*t^18 +12*z^11*t^20 -6*z^12*t^22 -2*z^2*t^4 -189*z^5*t^8
   +12*z^4*t^7 +37*z^3*t^5 -7*z^2*t^2 -z^8*t^16 +5*z^10*t^19 -6*z^5*t^10 -67*z^5*t^9
   -156*z^9*t^15 +39*z^10*t^17 +28*z^11*t^19 -8*z^12*t^21 +75*t^10*z^7 +t^12*z^8
   -85*t^14*z^9 -2*z^12*t^20 +12*t^18*z^11 -104*z^5*t^7 +17*z^10*t^16 -4*z^6*t^8
   -12*t^6*z^5 +2*z^4*t^4)/(1 -592*z^9*t^16 -105*z^6*t^10 -175*z^8*t^14 +8*z^7*t^14
   +192*t^6*z^4 +4*z^6*t^11 +935*t^12*z^7 -33*z^8*t^10 -65*z^9*t^17 -118*z^9*t^12
   +2*z^6*t^12 -202*z^8*t^13 +50*z^7*t^8 +154*z^7*t^13 -44*z^8*t^15 +2152*z^7*t^11
   +207*z^4*t^5 -426*t^9*z^6 -23*z^2*t^3 +17*t^3*z^3 +8*z^3*t^6 -z -2*z*t^2 -3*z*t
   +88*z^3*t^4 -88*z^6*t^7 -16*z^12*t^18 +2*z^13*t^20 +110*t^14*z^10 +7*t^16*z^11
   +10*z^11*t^21 -5*z^12*t^23 +203*z^10*t^18 +157*z^11*t^20 -69*z^12*t^22 -2*z^9*t^18
   -3*z^2*t^4 -607*z^5*t^8 +41*z^4*t^7 +2*z^4*t^8 +60*z^3*t^5 -20*z^2*t^2 -2*z^3*t^2
   -4*z^2*t -3*z^8*t^16 +27*z^10*t^19 +z^10*t^20 -12*z^5*t^10 -156*z^5*t^9
   -1992*z^9*t^15 +542*z^10*t^17 +755*z^11*t^19 -273*z^12*t^21 -12*z^13*t^24
   +6*z^14*t^26 +4*z^4*t^3 +z^5*t^4 +44*t^25*z^14 -100*t^23*z^13 +1612*t^10*z^7
   -41*t^12*z^8 -2484*t^14*z^9 -329*z^12*t^20 +54*z^14*t^24 +1208*t^18*z^11 -720*z^5*t^7
   -194*z^13*t^22 +611*z^10*t^16 -219*z^6*t^8 -284*t^6*z^5 -206*t^11*z^8 +673*z^11*t^17
   -88*z^13*t^21 +268*t^15*z^10 +16*t^23*z^14 -178*z^12*t^19 -922*z^9*t^13 +519*t^9*z^7
   +68*z^4*t^4 -13*z^6*t^6 -22*z^5*t^5)
\end{verbatim}
and
\begin{verbatim}
T_8(2,z,1) = (1 -16*z^12 +74*z^4 +597*z^7 -114*z^8 -336*z^9 -31*z^6 +52*z^11 +94*z^10
   +88*z^3 -22*z^2 -378*z^5 -5*z)/(1 +120*z^14 -870*z^12 +514*z^4 +5430*z^7 -704*z^8
   -6175*z^9 -845*z^6 +2810*z^11 +1762*z^10 +171*z^3 -50*z^2 -1800*z^5 -6*z -392*z^13)
\end{verbatim}

A subset of these results where $n=m$ collects the
of ways of placing $2\times 2$ squares into other squares \cite[A193580,A063443]{EIS}:

\small
\begin{Verbatim}
 2 1 1:  1 : 1
 2 2 2:  1 1 : 2
 2 3 3:  1 4 0 : 5
 2 4 4:  1 9 16 8 1 : 35
 2 5 5:  1 16 78 140 79 0 0 : 314
 2 6 6:  1 25 228 964 1987 1974 978 242 27 1 : 6427
 2 7 7:  1 36 520 3920 16834 42368 62266 51504 21792 3600 0 0 0 : 202841
 2 8 8:  1 49 1020 11860 85275 397014 1220298 2484382 3324193 2882737 \
          1601292 569818 129657 18389 1520 64 1 : 12727570
\end{Verbatim}
\normalsize

\subsection{$3\times 3$ Squares}
Placing $3\times 3$ squares into rectangles yields:
\small
\begin{Verbatim}
 3 1 1:  1 : 1
 3 2 1:  1 : 1
 3 2 2:  1 : 1
 3 3 1:  1 : 1
 3 3 2:  1 : 1
 3 3 3:  1 1 : 2
 3 4 1:  1 : 1
 3 4 2:  1 : 1
 3 4 3:  1 2 : 3
 3 4 4:  1 4 : 5
 3 5 1:  1 : 1
 3 5 2:  1 0 : 1
 3 5 3:  1 3 : 4
 3 5 4:  1 6 0 : 7
 3 5 5:  1 9 0 : 10
 3 6 1:  1 : 1
 3 6 2:  1 0 : 1
 3 6 3:  1 4 1 : 6
 3 6 4:  1 8 4 : 13
 3 6 5:  1 12 9 0 : 22
 3 6 6:  1 16 30 12 1 : 60
\end{Verbatim}
\normalsize
If the rectangle is only 3 units wide, the problem
is equivalent to tiling a $1\times m$ board with monomers and straight trimers,
see (\ref{eq.sonce}) and \cite[A102547]{EIS}:
\small
\begin{Verbatim}
 3 3 1:  1 : 1
 3 3 2:  1 : 1
 3 3 3:  1 1 : 2
 3 3 4:  1 2 : 3
 3 3 5:  1 3 : 4
 3 3 6:  1 4 1 : 6
 3 3 7:  1 5 3 : 9
 3 3 8:  1 6 6 : 13
 3 3 9:  1 7 10 1 : 19
 3 3 10:  1 8 15 4 : 28
 3 3 11:  1 9 21 10 : 41
 3 3 12:  1 10 28 20 1 : 60
 3 3 13:  1 11 36 35 5 : 88
 3 3 14:  1 12 45 56 15 : 129
\end{Verbatim}
\normalsize
The row sums are \cite[A000930]{EIS}
\begin{eqnarray}
T_{3\times m}(3) &=& T_{3\times (m-1)}(3)+ T_{3\times (m-3)}(3).\\
T_{s\times m}(s) &=& T_{s\times (m-1)}(s)+ T_{s\times (m-3)}(s). \label{eq.onelane}
\end{eqnarray}
The generating function is
\begin{eqnarray}
T_3(3,z,t) &=& \frac{1}{1-z-z^3t};\\
T_s(s,z,t) &=& \frac{1}{1-z-z^st}.
\label{eq.gss}
\end{eqnarray}
If the rectangle is 4 or 5 units wide, each square has one or two more 
places to go:
$T_{4\times m}(3,k)= 2^k T_{3\times m}(3,k)$
with generating functions
\begin{equation}
T_4(3,z,t) = \frac{1}{1-z-2z^3t};\quad 
T_5(3,z,t) = \frac{1}{1-z-3z^3t}.
\end{equation}
More generally one may account for the additional freedom with a
factor 
$n-s+1$ 
for each of the $k$ squares
if the width remains smaller than twice the square's size:
\begin{equation}
T_{n\times m}(s,k)= (n-s+1)^k T_{s\times m}(s,k), \quad s\le n<2s.
\label{eq.subwidth}
\end{equation}
This is echoed in the generating function (\ref{eq.2g}):
\begin{equation}
T_n(s,z,t)= \frac{1}{1-z-(n-s+1)z^st}, \quad s\le n<2s.
\label{eq.gfsub}
\end{equation}
If the rectangle is at least
twice as wide as the square, $n=2s$, squares may be
stacked along the short direction:
\small
\begin{Verbatim}
 3 6 1:  1 : 1
 3 6 2:  1 0 : 1
 3 6 3:  1 4 1 : 6
 3 6 4:  1 8 4 : 13
 3 6 5:  1 12 9 0 : 22
 3 6 6:  1 16 30 12 1 : 60
 3 6 7:  1 20 67 50 9 : 147
 3 6 8:  1 24 120 128 36 0 : 309
 3 6 9:  1 28 189 310 166 26 1 : 721
 3 6 10:  1 32 274 660 561 176 16 : 1720
 3 6 11:  1 36 375 1242 1461 672 100 0 : 3887
 3 6 12:  1 40 492 2120 3362 2236 600 48 1 : 8900
 3 6 13:  1 44 625 3358 7016 6480 2721 470 25 : 20740
 3 6 14:  1 48 774 5020 13431 16296 9438 2472 225 0 : 47705
 3 6 15:  1 52 939 7170 23871 36880 28220 10582 1713 80 1 : 109509
\end{Verbatim}
\normalsize
with generating function
\begin{verbatim}
T_6(3,z,t) = ( -z^2*t +1 -z^3*t^2)/( -3*z^3*t -z^4*t^2 +1 -z^2*t -2*z^3*t^2 
                              -z +z^5*t^3 +z^6*t^4 +2*z^5*t^2 +2*z^6*t^3),
\end{verbatim}
and the associated generating function of the row sums
\begin{verbatim}
T_6(3,z,1) = ( -z^2 +1 -z^3)/( -5*z^3 -z^4 +1 -z^2 -z +3*z^5 +3*z^6).
\end{verbatim}
If the rectangle is 7 units wide
\small
\begin{Verbatim}
 3 7 1:  1 : 1
 3 7 2:  1 0 : 1
 3 7 3:  1 5 3 : 9
 3 7 4:  1 10 12 0 : 23
 3 7 5:  1 15 27 0 : 43
 3 7 6:  1 20 67 50 9 : 147
 3 7 7:  1 25 132 200 79 0 : 437
 3 7 8:  1 30 222 500 314 0 0 : 1067
 3 7 9:  1 35 337 1075 1179 333 27 0 : 2987
 3 7 10:  1 40 477 2050 3469 2160 408 0 : 8605
 3 7 11:  1 45 642 3550 8309 7998 2508 0 0 : 23053
 3 7 12:  1 50 832 5700 17449 23936 13018 1820 81 0 : 62887
 3 7 13:  1 55 1047 8625 33264 61599 52089 17218 1847 0 0 : 175745
 3 7 14:  1 60 1287 12450 58754 140112 165540 87852 16147 0 0 : 482203
\end{Verbatim}
\normalsize
with generating function
\begin{verbatim}
T_7(3,z,t) = (z^6*t^4 -z^5*t^3 -3*z^3*t^2 -z^4*t^2 -z^2*t +1)/( -3*z^9*t^6 
    +3*z^8*t^5 +4*z^7*t^4 +10*z^6*t^4 +z^7*t^3 +z^8*t^4 +2*z^5*t^3 -4*z^4*t^2 
    -6*z^3*t^2 +2*z^6*t^3 -4*z^3*t -z^2*t -z^9*t^5 -z +1)
\end{verbatim}
and generating function
\begin{verbatim}
T_7(3,z,1) = (z^6 -z^5 -3*z^3 -z^4 -z^2 +1)/( -4*z^9 +4*z^8 +5*z^7 
                              +12*z^6 +2*z^5 -4*z^4 -10*z^3 -z^2 -z +1)
\end{verbatim}
for the row sums.
If the rectangle is 8 units wide
\small
\begin{Verbatim}
 3 8 1:  1 : 1
 3 8 2:  1 0 : 1
 3 8 3:  1 6 6 : 13
 3 8 4:  1 12 24 0 : 37
 3 8 5:  1 18 54 0 0 : 73
 3 8 6:  1 24 120 128 36 0 : 309
 3 8 7:  1 30 222 500 314 0 0 : 1067
 3 8 8:  1 36 360 1232 1246 0 0 0 : 2875
\end{Verbatim}
\normalsize
with generating function
\begin{verbatim}
T_8(3,z,t) = ( -z^8*t^5 +2*z^7*t^4 +z^5*t^3 -6*z^3*t^2 +1 +5*z^6*t^4 -2*z^2*t 
    -z^4*t^2 +z^6*t^3 -z^9*t^6)/(6*z^12*t^8 +6*z^11*t^7 -13*z^10*t^6 -31*z^9*t^6 
    -7*z^9*t^5 -7*z^8*t^5 +2*z^8*t^4 +13*z^7*t^4 +41*z^6*t^4 +t^3*z^7 +4*z^6*t^3 
    +13*z^5*t^3 +z^5*t^2 -7*z^4*t^2 -12*z^3*t^2 -4*z^3*t -2*z^2*t +1 -z)
\end{verbatim}
with row sums
\begin{verbatim}
T_8(3,z,1) = ( -z^8 +2*z^7 +z^5 -6*z^3 +1 +6*z^6 -2*z^2 -z^4 -z^9) /(6*z^12 
    +6*z^11 -13*z^10 -38*z^9 -5*z^8 +14*z^7 +45*z^6 +14*z^5 -7*z^4 -16*z^3 -2*z^2 +1 -z).
\end{verbatim}
If the rectangle is 9 units wide
\small
\begin{Verbatim}
3 9 0: 1 : 1
3 9 1: 1 0 : 1
3 9 2: 1 0 0 : 1
3 9 3: 1 7 10 1 : 19
3 9 4: 1 14 40 8 0 : 63
3 9 5: 1 21 90 27 0 0 : 139
3 9 6: 1 28 189 310 166 26 1 : 721
3 9 7: 1 35 337 1075 1179 333 27 0 : 2987
3 9 8: 1 42 534 2540 4316 1740 216 0 0 : 9389
3 9 9: 1 49 780 5048 13211 11984 4526 758 51 1 : 36409
3 9 10: 1 56 1075 8942 33356 53062 37007 11116 1444 64 0 : 146123
\end{Verbatim}
\normalsize
with generating function
\begin{verbatim}
T_9(3,z,t) = (1 +2*z^8*t^4 -6*z^10*t^6 -12*z^3*t^2 -2*z^2*t +5*t^10*z^14 -5*z^3*t^3
   -2*z^5*t^3 +10*z^6*t^6 -t^14*z^18 +40*z^6*t^5 -z^15*t^15 -19*z^11*t^9
   +32*z^12*t^10 -50*z^10*t^8 -53*z^9*t^7 -2*z^9*t^6 +23*z^13*t^10 +6*z^14*t^11
   +13*t^5*z^8 +12*t^4*z^7 -4*z^12*t^9 -32*t^8*z^11 +z^17*t^16 -z^21*t^19 +3*z^18*t^17
   +z^19*t^17 -z^20*t^17 +7*z^18*t^16 +z^20*t^18 +5*z^12*t^12 -19*z^15*t^14
   -57*z^9*t^8 -29*t^7*z^10 +12*z^7*t^6 -10*z^9*t^9 +12*z^13*t^11 +13*z^14*t^12
   -3*z^10*t^9 +8*z^11*t^10 +45*z^12*t^11 -5*z^14*t^13 -2*z^13*t^12 -4*z^4*t^2
   +40*z^7*t^5 +8*z^8*t^6 +32*z^6*t^4 +4*z^9*t^5 -7*z^4*t^3 +5*z^5*t^4 -7*z^8*t^7
   +2*t^3*z^6 +5*t^14*z^17 +t^13*z^17 -14*t^13*z^15 +t^14*z^16 +t^15*z^18 -4*t^13*z^16
   +2*t^11*z^15 +2*t^12*z^15 -4*z^17*t^15 -4*z^13*t^9 -4*z^11*t^7 -t^8*z^12)/(1
   +14*z^8*t^4 -110*z^10*t^6 -22*z^3*t^2 +34*t^9*z^14 -2*z^2*t +264*t^10*z^14
   -6*z^3*t^3 -5*z^3*t +52*t^8*z^13 +29*t^11*z^16 +13*z^5*t^3 +15*z^6*t^6
   -19*t^14*z^18 +6*t^10*z^15 -z +2*z^9*t^4 -3*z^21*t^20 -z^23*t^21 -6*z^23*t^20
   -6*z^22*t^19 +7*z^24*t^21 +7*z^21*t^16 -5*z^21*t^17 +7*z^23*t^19 +3*z^16*t^15
   -7*z^20*t^15 +96*z^6*t^5 -6*z^15*t^15 +34*z^11*t^9 +526*z^12*t^10 -238*z^10*t^8
   +z^18*t^18 -38*t^12*z^17 -435*z^9*t^7 -2*z^5*t^2 -282*z^9*t^6 +446*z^13*t^10
   +190*z^14*t^11 +70*t^5*z^8 +84*t^4*z^7 +417*z^12*t^9 -120*t^8*z^11 -36*z^20*t^16
   -z^22*t^20 -51*z^21*t^18 +6*z^17*t^16 -29*z^21*t^19 +29*z^18*t^17 -3*z^19*t^17
   +21*z^20*t^17 +153*z^18*t^16 -2*z^20*t^18 +15*z^12*t^12 -99*z^15*t^14 -z^20*t^19
   +z^24*t^22 -15*t^13*z^18 -176*z^9*t^8 -431*t^7*z^10 +z^19*t^14 +27*z^7*t^6
   -20*z^9*t^9 +114*z^13*t^11 -25*z^14*t^12 -23*z^10*t^9 +15*z^11*t^10 +175*z^12*t^11
   -13*z^14*t^13 +2*z^13*t^12 -12*z^4*t^2 +145*z^7*t^5 -28*z^8*t^6 +142*z^6*t^4
   -38*z^12*t^7 -23*z^11*t^6 +4*z^7*t^3 -z^9*t^5 -8*z^4*t^3 +7*z^5*t^4 -12*z^8*t^7
   +4*t^3*z^6 -102*t^14*z^17 -165*t^12*z^16 -42*t^13*z^17 -371*t^13*z^15 +4*t^14*z^16
   +98*t^15*z^18 -4*t^16*z^19 +25*t^15*z^19 -122*t^13*z^16 +14*t^11*z^15 -211*t^12*z^15
   +18*z^17*t^15 +254*z^13*t^9 -147*z^11*t^7 -t^8*z^12 +4*t^5*z^10)
\end{verbatim}
with row sums
\begin{verbatim}
T_9(3,z,1) = (10*z^18 -30*z^15 -47*z^11 +19*z^14 +3*z^5 +84*z^6 +16*z^8 -88*z^10
   -17*z^3 -2*z^2 -z^21 +1 +z^19 +3*z^17 +29*z^13 +64*z^7 +77*z^12 -118*z^9 -11*z^4
   -3*z^16) /(247*z^18 -667*z^15 -241*z^11 +450*z^14 +18*z^5 +257*z^6 +44*z^8 -798*z^10
   -33*z^3 -2*z^2 -81*z^21 +8*z^24 +1 -z +19*z^19 -25*z^20 -158*z^17 +868*z^13 +260*z^7
   +1094*z^12 -912*z^9 -20*z^4 -251*z^16 -7*z^22)
\end{verbatim}

To summarize, here is the number
of ways of placing $3\times 3$ squares into other squares \cite[A276171]{EIS}:
\small
\begin{Verbatim}
 3 1 1:  1 : 1
 3 2 2:  1 : 1
 3 3 3:  1 1 : 2
 3 4 4:  1 4 : 5
 3 5 5:  1 9 0 : 10
 3 6 6:  1 16 30 12 1 : 60
 3 7 7:  1 25 132 200 79 0 : 437
 3 8 8:  1 36 360 1232 1246 0 0 0 : 2875
 3 9 9:  1 49 780 5048 13211 11984 4526 758 51 1 : 36409
 3 10 10:  1 64 1470 15468 78851 193672 234394 139188 37760 3600 0 0 : 704468
 3 11 11:  1 81 2520 38972 324721 1490562 3761236 5052890 3305328 807648 : 14783959
 3 12 12:  1 100 4032 85600 1050442 7728696 34400276 91562420 141084672 \
           119132920 52175594 12725724 1828210 152908 6884 144 1 : 461938624
\end{Verbatim}
\normalsize
More row sums are in \cite[A140304]{EIS}.

The geometries of placing squares into squares of twice
the edge length, $n=m=2s$, are with (\ref{eq.k1})
\begin{equation}
T_{2s\times 2s}(s,1)=(s+1)^2,
\end{equation}
with (\ref{eq.ij})
\begin{equation}
T_{2s\times 2s}(s,4)=1,
\end{equation}
and otherwise counted by considering the few
number of constellations where all squares touch the
bigger square:
\begin{equation}
T_{2s\times 2s}(s,2)=2s(s+2), \quad
T_{2s\times 2s}(s,3)=4s.
\end{equation}

\subsection{$4\times 4$ Squares}
Placing $4\times 4$ squares into rectangles yields:
\small
\begin{Verbatim}
 4 1 1:  1 : 1
 4 2 1:  1 : 1
 4 2 2:  1 : 1
 4 3 1:  1 : 1
 4 3 2:  1 : 1
 4 3 3:  1 : 1
 4 4 1:  1 : 1
 4 4 2:  1 : 1
 4 4 3:  1 : 1
 4 4 4:  1 1 : 2
 4 5 1:  1 : 1
 4 5 2:  1 : 1
 4 5 3:  1 : 1
 4 5 4:  1 2 : 3
 4 5 5:  1 4 : 5
 4 6 1:  1 : 1
 4 6 2:  1 : 1
 4 6 3:  1 0 : 1
 4 6 4:  1 3 : 4
 4 6 5:  1 6 : 7
 4 6 6:  1 9 0 : 10
 4 7 1:  1 : 1
 4 7 2:  1 : 1
 4 7 3:  1 0 : 1
 4 7 4:  1 4 : 5
 4 7 5:  1 8 0 : 9
 4 7 6:  1 12 0 : 13
 4 7 7:  1 16 0 0 : 17
\end{Verbatim}
\normalsize
If the rectangle is 4 units wide, the problem is equivalent
to placing 1-ominos and straight tetrominos on  a line, see (\ref{eq.sonce}) and  \cite[A180184]{EIS}:
\small
\begin{Verbatim}
 4 4 1:  1 : 1
 4 4 2:  1 : 1
 4 4 3:  1 : 1
 4 4 4:  1 1 : 2
 4 4 5:  1 2 : 3
 4 4 6:  1 3 : 4
 4 4 7:  1 4 : 5
 4 4 8:  1 5 1 : 7
 4 4 9:  1 6 3 : 10
 4 4 10:  1 7 6 : 14
 4 4 11:  1 8 10 : 19
 4 4 12:  1 9 15 1 : 26
 4 4 13:  1 10 21 4 : 36
 4 4 14:  1 11 28 10 : 50
 4 4 15:  1 12 36 20 : 69
 4 4 16:  1 13 45 35 1 : 95
 4 4 17:  1 14 55 56 5 : 131
 4 4 18:  1 15 66 84 15 : 181
 4 4 19:  1 16 78 120 35 : 250
\end{Verbatim}
\normalsize
The generating function is
\begin{equation}
T_4(4,z,t) = \frac{1}{1-z-z^4t}.
\end{equation}
$T_n(4,z,t)$ in the range $n=5\ldots 7$ are given by increasing
the factor in front of the $t$ in the denominator,
see (\ref{eq.subwidth}) and (\ref{eq.gfsub}).
The row sums are given by (\ref{eq.onelane}) and \cite[A003269]{EIS}.
If the rectangle is 8 units wide,
\small
\begin{Verbatim}
 4 8 1:  1 : 1
 4 8 2:  1 0 : 1
 4 8 3:  1 0 : 1
 4 8 4:  1 5 1 : 7
 4 8 5:  1 10 4 : 15
 4 8 6:  1 15 9 0 : 25
 4 8 7:  1 20 16 0 : 37
 4 8 8:  1 25 48 16 1 : 91
 4 8 9:  1 30 105 66 9 : 211
 4 8 10:  1 35 187 168 36 0 : 427
 4 8 11:  1 40 294 340 100 0 : 775
 4 8 12:  1 45 426 707 342 39 1 : 1561
 4 8 13:  1 50 583 1394 1031 252 16 : 3327
 4 8 14:  1 55 765 2526 2564 938 100 0 : 6949
 4 8 15:  1 60 972 4228 5466 2628 400 0 : 13755
 4 8 16:  1 65 1204 6625 10887 6905 1621 82 1 : 27391
 4 8 17:  1 70 1461 9842 20602 17220 6002 740 25 : 55963
 4 8 18:  1 75 1743 14004 37011 39708 18876 3708 225 0 : 115351
 4 8 19:  1 80 2050 19236 63139 84004 50814 13548 1225 0 : 234097
 4 8 20:  1 85 2382 25663 102636 166368 124806 43544 5992 155 1 : 471633
 4 8 21:  1 90 2739 33410 159777 312810 288340 129780 26346 1842 36 : 955171
 4 8 22:  1 95 3121 42602 239462 562215 630379 357750 99702 11823 441 0 : 1947591
\end{Verbatim}
\normalsize
with generating function
\begin{verbatim}
T_8(4,z,t) = ( -z^2*t +1 +z^6*t^3 -z^4*t^2 -z^3*t)/(2*z^7*t^3 +3*z^7*t^2 -z 
    -z^5*t^2 -4*z^4*t -z^2*t +1 +2*z^6*t^3 -2*z^4*t^2 -z^10*t^5 +z^8*t^4 +3*z^6*t^2 
    -3*z^10*t^4 +3*z^8*t^3)
\end{verbatim}
and associated generating function of the row sums
\begin{verbatim}
T_8(4,z,1) = ( -z^2 +1 +z^6 -z^4 -z^3)
    /(5*z^7 -z -z^5 -6*z^4 -z^2 +1 +5*z^6 -4*z^10 +4*z^8).
\end{verbatim}

If the rectangle is 9 units wide,
\begin{verbatim}
T_9(4,z,t) = ( -z^7*t^3 -4*z^4*t^2 -z^6*t^2 -2*z^5*t^2 +4*z^8*t^4 +3*z^9*t^4 
    -z^10*t^5 +3*z^6*t^3 -z^3*t -z^2*t +1 -z^12*t^6)/(1 +16*z^8*t^4 +7*z^8*t^3 
    +6*z^6*t^3 -z^2*t +5*z^7*t^3 +2*z^11*t^5 +z^6*t^2 -z +13*z^9*t^4 -7*z^4*t^2 
    -4*z^5*t^2 -5*z^4*t +4*t^3*z^9 +z^7*t^2 -13*z^12*t^6 -10*z^10*t^5 -8*z^12*t^5 
    +2*z^14*t^6 -10*z^13*t^6 +3*z^14*t^7 +3*z^16*t^8 -6*z^13*t^5 +2*z^16*t^7 
    +2*z^11*t^4 +2*z^10*t^3)
\end{verbatim}
with row sums
\begin{verbatim}
T_9(4,z,1) = ( -z^7 -4*z^4 +2*z^6 -2*z^5 +4*z^8 +3*z^9 -z^10 -z^3 -z^2 +1 -z^12)
    /(1 -z -z^2 -12*z^4 -4*z^5 +7*z^6 +23*z^8 +17*z^9 +6*z^7 +4*z^11 -8*z^10 
    -21*z^12 +5*z^14 -16*z^13 +5*z^16).
\end{verbatim}

In overview, this is the number
of ways of placing $4\times 4$ squares into other squares:
\small
\begin{Verbatim}
 4 1 1:  1 : 1
 4 2 2:  1 : 1
 4 3 3:  1 : 1
 4 4 4:  1 1 : 2
 4 5 5:  1 4 : 5
 4 6 6:  1 9 0 : 10
 4 7 7:  1 16 0 0 : 17
 4 8 8:  1 25 48 16 1 : 91
 4 9 9:  1 36 198 260 79 0 : 574
 4 10 10:  1 49 516 1568 1246 0 0 : 3380
 4 11 11:  1 64 1080 6000 9550 0 0 0 : 16695
 4 12 12:  1 81 1980 18200 59974 48408 14412 1764 81 1 : 144902
 4 13 13:  1 100 3318 46084 273197 637936 644318 292192 56560 3600 0 : 1957306
 4 14 14:  1 121 5208 101860 957384 4296934 9399764 9966546 4789008 807648 0 0 0 : 30324474
\end{Verbatim}
\normalsize

\subsection{$5\times 5$ Squares}
Placing $5\times 5$ squares into rectangles yields:
\small
\begin{Verbatim}
 5 1 1:  1 : 1
 5 2 1:  1 : 1
 5 2 2:  1 : 1
 5 3 1:  1 : 1
 5 3 2:  1 : 1
 5 3 3:  1 : 1
 5 4 1:  1 : 1
 5 4 2:  1 : 1
 5 4 3:  1 : 1
 5 4 4:  1 : 1
 5 5 1:  1 : 1
 5 5 2:  1 : 1
 5 5 3:  1 : 1
 5 5 4:  1 : 1
 5 5 5:  1 1 : 2
 5 6 1:  1 : 1
 5 6 2:  1 : 1
 5 6 3:  1 : 1
 5 6 4:  1 : 1
 5 6 5:  1 2 : 3
 5 6 6:  1 4 : 5
 5 7 1:  1 : 1
 5 7 2:  1 : 1
 5 7 3:  1 : 1
 5 7 4:  1 0 : 1
 5 7 5:  1 3 : 4
 5 7 6:  1 6 : 7
 5 7 7:  1 9 : 10
 5 8 1:  1 : 1
 5 8 2:  1 : 1
 5 8 3:  1 : 1
 5 8 4:  1 0 : 1
 5 8 5:  1 4 : 5
 5 8 6:  1 8 : 9
 5 8 7:  1 12 0 : 13
 5 8 8:  1 16 0 : 17
 5 9 1:  1 : 1
 5 9 2:  1 : 1
 5 9 3:  1 0 : 1
 5 9 4:  1 0 : 1
 5 9 5:  1 5 : 6
 5 9 6:  1 10 0 : 11
 5 9 7:  1 15 0 : 16
 5 9 8:  1 20 0 : 21
 5 9 9:  1 25 0 0 : 26
\end{Verbatim}
\normalsize
If the rectangle is 5 units wide, the problem is equivalent
to placing 1-ominos and straight pentominos on  a line, 
see (\ref{eq.sonce}), (\ref{eq.gss}) and \cite[A003520]{EIS}:
\small
\begin{Verbatim}
 5 5 1:  1 : 1
 5 5 2:  1 : 1
 5 5 3:  1 : 1
 5 5 4:  1 : 1
 5 5 5:  1 1 : 2
 5 5 6:  1 2 : 3
 5 5 7:  1 3 : 4
 5 5 8:  1 4 : 5
 5 5 9:  1 5 : 6
 5 5 10:  1 6 1 : 8
 5 5 11:  1 7 3 : 11
 5 5 12:  1 8 6 : 15
 5 5 13:  1 9 10 : 20
 5 5 14:  1 10 15 : 26
 5 5 15:  1 11 21 1 : 34
 5 5 16:  1 12 28 4 : 45
 5 5 17:  1 13 36 10 : 60
 5 5 18:  1 14 45 20 : 80
 5 5 19:  1 15 55 35 : 106
\end{Verbatim}
\normalsize
If the rectangle is 6 to 9 units wide, the counts
are described by (\ref{eq.subwidth}) and (\ref{eq.gfsub}). 
If it is 10 units wide,
\small
\begin{Verbatim}
5 10 1: 1 : 1
5 10 2: 1 : 1
5 10 3: 1 0 : 1
5 10 4: 1 0 : 1
5 10 5: 1 6 1 : 8
5 10 6: 1 12 4 : 17
5 10 7: 1 18 9 : 28
5 10 8: 1 24 16 0 : 41
5 10 9: 1 30 25 0 : 56
5 10 10: 1 36 70 20 1 : 128
5 10 11: 1 42 151 82 9 : 285
5 10 12: 1 48 268 208 36 : 561
\end{Verbatim}
\normalsize
with generating function
\begin{verbatim}
T_10(5,z,t) = ( -2*z^5*t^2 -z^3*t +1 +z^10*t^4 +z^8*t^3 -z^6*t^2 -z^4*t)
    /( -z^15*t^6 -4*z^15*t^5 -t^5*z^13 -4*z^13*t^4 +2*z^11*t^4 +3*z^10*t^4
     +4*z^11*t^3 +8*z^10*t^3 +2*z^9*t^3 +2*z^8*t^3 +4*z^9*t^2 +4*z^8*t^2
     -z^7*t^2 -z^6*t^2 -3*z^5*t^2 -5*z^5*t -z^3*t +1 -z)
\end{verbatim}
and for row sums
\begin{verbatim}
T_10(5,z,1) = ( -2*z^5 -z^3 +1 +z^10 +z^8 -z^6 -z^4)/( -5*z^15 -5*z^13 +6*z^11 
    +11*z^10 +6*z^9 +6*z^8 -z^7 -z^6 -8*z^5 -z^3 +1 -z).
\end{verbatim}
If the rectangle is 11 units wide,
\small
\begin{Verbatim}
5 11 1: 1 : 1
5 11 2: 1 : 1
5 11 3: 1 0 : 1
5 11 4: 1 0 : 1
5 11 5: 1 7 3 : 11
5 11 6: 1 14 12 : 27
5 11 7: 1 21 27 0 : 49
5 11 8: 1 28 48 0 : 77
5 11 9: 1 35 75 0 : 111
5 11 10: 1 42 151 82 9 : 285
5 11 11: 1 49 276 320 79 : 725
5 11 12: 1 56 450 788 314 0 : 1609
\end{Verbatim}
\normalsize
\begin{verbatim}
with generating function
T_11(5,z,t) = (11*z^10*t^4 +2*z^8*t^3 -z^8*t^2 -z^10*t^3 +z^20*t^8 +4*z^12*t^4 
    -6*z^15*t^6 -3*z^9*t^3 -z^4*t +2*z^13*t^5 -4*z^6*t^2 -4*z^16*t^6 +7*z^11*t^4 
    -2*z^7*t^2 +1 -6*z^5*t^2 -z^3*t -z^18*t^7 +3*z^14*t^5)/(1 -z -4*z^13*t^5
    -6*z^5*t -9*z^5*t^2 -4*z^7*t^2 +16*z^10*t^3 +14*z^14*t^5 +28*z^11*t^4
    -4*z^6*t^2 +29*z^10*t^4 +4*z^9*t^3 -39*z^15*t^6 -z^3*t +5*z^8*t^3 -27*z^15*t^5
    +3*z^15*t^4 +9*z^14*t^4 -31*z^16*t^6 +19*z^20*t^8 +6*z^12*t^3 +11*z^11*t^3
    -21*z^16*t^5 +3*z^23*t^9 +18*z^20*t^7 -10*z^19*t^7 -7*z^18*t^7 -9*z^19*t^6
    -6*z^18*t^6 +3*z^23*t^8 +z^13*t^4 -3*z^25*t^10 -16*z^17*t^6 +13*z^21*t^8
    +2*z^8*t^2 +15*z^12*t^4 -3*z^25*t^9 +12*z^21*t^7 +2*z^9*t^2 +3*z^13*t^3
    -12*z^17*t^5)
\end{verbatim}
and with row sums
\begin{verbatim}
T_11(5,z,1) = (10*z^10 +z^8 +z^20 +4*z^12 -6*z^15 -3*z^9 -z^4 +2*z^13 -4*z^6
    -4*z^16 +7*z^11 -2*z^7 +1 -6*z^5 -z^3 -z^18 +3*z^14)/(1 +37*z^20 -6*z^25 -z
    +6*z^23 -19*z^19 -13*z^18 +23*z^14 +7*z^8 +21*z^12 -63*z^15 +6*z^9 -4*z^6
    -52*z^16 +45*z^10 +39*z^11 -4*z^7 -15*z^5 -z^3 -28*z^17 +25*z^21).
\end{verbatim}
Here is the number
of ways of placing $5\times 5$ squares into other squares:

\small
\begin{Verbatim}
 5 1 1:  1 : 1
 5 2 2:  1 : 1
 5 3 3:  1 : 1
 5 4 4:  1 : 1
 5 5 5:  1 1 : 2
 5 6 6:  1 4 : 5
 5 7 7:  1 9 : 10
 5 8 8:  1 16 0 : 17
 5 9 9:  1 25 0 0 : 26
 5 10 10:  1 36 70 20 1 : 128
 5 11 11:  1 49 276 320 79 : 725
 5 12 12:  1 64 696 1904 1246 0 : 3911
 5 13 13:  1 81 1420 7200 9550 0 0 : 18252
 5 14 14:  1 100 2550 20900 48175 0 0 0 : 71726
 5 15 15:  1 121 4200 52140 214680 153190 37040 3476 117 1 : 464966
 5 16 16:  1 144 6496 115104 790396 1729976 1479306 532572 78192 3600 0 : 4735787
\end{Verbatim}
\normalsize

\subsection{$6\times 6$ Squares}
If the rectangle is 6 units wide, the problem is equivalent
to placing 1-ominos and straight 6-ominos on  a line, see (\ref{eq.sonce}) and \cite[A005708]{EIS}:
\small
\begin{Verbatim}
 6 6 1:  1 : 1
 6 6 2:  1 : 1
 6 6 3:  1 : 1
 6 6 4:  1 : 1
 6 6 5:  1 : 1
 6 6 6:  1 1 : 2
 6 6 7:  1 2 : 3
 6 6 8:  1 3 : 4
 6 6 9:  1 4 : 5
 6 6 10:  1 5 : 6
 6 6 11:  1 6 : 7
 6 6 12:  1 7 1 : 9
 6 6 13:  1 8 3 : 12
 6 6 14:  1 9 6 : 16
 6 6 15:  1 10 10 : 21
 6 6 16:  1 11 15 : 27
 6 6 17:  1 12 21 : 34
 6 6 18:  1 13 28 1 : 43
\end{Verbatim}
\normalsize
If the rectangle is 7 to 11 units wide, the counts
are described by (\ref{eq.subwidth}) and (\ref{eq.gfsub}). 
If the rectangle is 12 units wide,
\small
\begin{Verbatim}
6 12 1: 1 : 1
6 12 2: 1 : 1
6 12 3: 1 0 : 1
6 12 4: 1 0 : 1
6 12 5: 1 0 : 1
6 12 6: 1 7 1 : 9
6 12 7: 1 14 4 : 19
6 12 8: 1 21 9 : 31
6 12 9: 1 28 16 0 : 45
6 12 10: 1 35 25 0 : 61
6 12 11: 1 42 36 0 : 79
6 12 12: 1 49 96 24 1 : 171
6 12 13: 1 56 205 98 9 : 369
6 12 14: 1 63 363 248 36 : 711
6 12 15: 1 70 570 500 100 0 : 1241
6 12 16: 1 77 826 880 225 0 : 2009
\end{Verbatim}
\normalsize
with generating function
\begin{verbatim}
T_12(6,z,t) = ( -2*z^6*t^2 -z^4*t +2*z^9*t^3 +1 -z^3*t +z^12*t^4 +2*z^10*t^3
    -z^15*t^5 -z^7*t^2 -z^5*t)/(1 +z^21*t^7 +3*z^9*t^3 -z^3*t +2*z^13*t^4 -z^8*t^2
    +10*z^12*t^3 -5*z^18*t^5 +5*z^13*t^3 +3*z^11*t^3 -z^18*t^6 +5*z^11*t^2
    +5*z^21*t^6 -3*z^15*t^5 -10*z^15*t^4 +5*z^9*t^2 -z -z^7*t^2 -3*z^6*t^2
    +5*z^10*t^2 -6*z^6*t -3*z^16*t^5 +3*z^10*t^3 +3*z^12*t^4 -10*z^16*t^4)
\end{verbatim}
and with row sums
\begin{verbatim}
T_12(6,z,1) = ( -2*z^6 -z^4 +2*z^9 +1 -z^3 +z^12 +2*z^10 -z^15 -z^7 -z^5)
    /(1 +8*z^10 +6*z^21 -13*z^16 -6*z^18 -z -z^3 +8*z^9 -z^8 +7*z^13 -9*z^6
      -13*z^15 +8*z^11 -z^7 +13*z^12).
\end{verbatim}
If the rectangle is 13 units wide,
\small
\begin{Verbatim}
6 13 0: 1 : 1
6 13 1: 1 : 1
6 13 2: 1 : 1
6 13 3: 1 0 : 1
6 13 4: 1 0 : 1
6 13 5: 1 0 : 1
6 13 6: 1 8 3 : 12
6 13 7: 1 16 12 : 29
6 13 8: 1 24 27 : 52
6 13 9: 1 32 48 0 : 81
6 13 10: 1 40 75 0 : 116
6 13 11: 1 48 108 0 : 157
6 13 12: 1 56 205 98 9 : 369
6 13 13: 1 64 366 380 79 : 890
6 13 14: 1 72 591 932 314 0 : 1910
6 13 15: 1 80 880 1840 870 0 : 3671
6 13 16: 1 88 1233 3190 1955 0 : 6467
6 13 17: 1 96 1650 5068 3829 0 0 : 10644
6 13 18: 1 104 2131 7986 7968 906 27 : 19123
6 13 19: 1 112 2676 12456 16820 5232 408 : 37705
\end{Verbatim}
\normalsize
with generating function
\begin{verbatim}
T_13(6,z,t) = (1 +6*z^21*t^7 +7*z^24*t^8 -z^27*t^9 -18*z^19*t^6 -z^30*t^10
   +2*z^22*t^7 +13*z^14*t^4 +6*z^25*t^8 +16*z^13*t^4 -2*z^16*t^5 +2*z^17*t^5
   -11*z^20*t^6 -17*z^18*t^6 -11*z^15*t^5 -z^12*t^3 +6*z^15*t^4 -t^2*z^10
   -2*t^2*z^9 -z^3*t +6*z^9*t^3 -z^5*t -7*z^6*t^2 -5*z^7*t^2 -z^4*t +17*z^12*t^4
   -3*z^11*t^3 -3*z^8*t^2 +5*z^10*t^3)/(1 +4*t^3*z^16 -z +39*z^21*t^7 +58*z^24*t^8
   -19*z^27*t^9 -81*z^19*t^6 -22*z^30*t^10 -24*z^27*t^8 -64*z^19*t^5 +68*z^24*t^7
   -28*z^30*t^9 +8*z^22*t^6 +12*z^14*t^3 +72*z^25*t^7 +14*z^22*t^7 +30*z^14*t^4
   +67*z^25*t^8 +40*z^13*t^4 -26*z^16*t^5 +11*z^17*t^5 -66*z^20*t^6 +15*z^21*t^6
   -68*z^18*t^6 -29*z^15*t^5 +26*z^12*t^3 -5*z^15*t^4 +3*t^2*z^10 +3*t^2*z^9 -7*z^6*t
   -z^3*t +12*t^4*z^17 -8*t^8*z^28 -4*z^23*t^7 -8*t^6*z^23 +9*z^9*t^3 -63*z^18*t^5
   +19*z^13*t^3 -10*z^6*t^2 +39*z^26*t^8 +44*z^26*t^7 -24*t^5*z^21 -52*t^5*z^20
   -4*z^7*t^2 +3*z^33*t^11 +4*z^33*t^10 -11*t^4*z^16 +4*z^18*t^4 +3*z^36*t^12
   +4*z^36*t^11 -7*z^28*t^9 +38*z^12*t^4 +7*z^11*t^3 +3*t^2*z^11 -19*z^31*t^10
   -24*z^31*t^9 +8*z^15*t^3 -4*z^8*t^2 +8*z^10*t^3)
\end{verbatim}
and with row sums
\begin{verbatim}
T_13(6,z,1) = (1 -3*z^11 -3*z^8 +6*z^21 +7*z^24 -z^27 -18*z^19 +2*z^22 +13*z^14
   +6*z^25 +16*z^13 -2*z^16 +2*z^17 -11*z^20 -17*z^18 -5*z^15 +16*z^12 +4*z^10 +4*z^9
   -z^3 -z^5 -7*z^6 -5*z^7 -z^4 -z^30)/(1 -z +10*z^11 -4*z^8 +30*z^21 +126*z^24
   -43*z^27 -145*z^19 +22*z^22 +42*z^14 +139*z^25 +59*z^13 -33*z^16 +23*z^17 -118*z^20
   -127*z^18 -26*z^15 +64*z^12 +11*z^10 +12*z^9 -z^3 -17*z^6 -4*z^7 -50*z^30 +7*z^36
   -15*z^28 -12*z^23 +83*z^26 +7*z^33 -43*z^31) .
\end{verbatim}

This is the number
of ways of placing $6\times 6$ squares into other squares:

\small
\begin{Verbatim}
 6 1 1:  1 : 1
 6 2 2:  1 : 1
 6 3 3:  1 : 1
 6 4 4:  1 : 1
 6 5 5:  1 : 1
 6 6 6:  1 1 : 2
 6 7 7:  1 4 : 5
 6 8 8:  1 9 : 10
 6 9 9:  1 16 0 : 17
 6 10 10:  1 25 0 : 26
 6 11 11:  1 36 0 0 : 37
 6 12 12:  1 49 96 24 1 : 171
 6 13 13:  1 64 366 380 79 : 890
 6 14 14:  1 81 900 2240 1246 0 : 4468
 6 15 15:  1 100 1800 8400 9550 0 0 : 19851
 6 16 16:  1 121 3180 24200 48175 0 0 0 : 75677
 6 17 17:  1 144 5166 58604 184681 0 0 0 0 : 248596
 6 18 18:  1 169 7896 127484 648301 410086 82594 6146 159 1 : 1282837
 6 19 19:  1 196 11520 253920 2016380 4105216 3013174 885104 102656 3600 0 : 10391767
\end{Verbatim}
\normalsize

\subsection{$7\times 7$ Squares}
If the rectangle is 7 units wide, the problem is equivalent
to placing 1-ominos and straight 7-ominos on  a line, see (\ref{eq.sonce}) and \cite[A005709]{EIS}:
\small
\begin{Verbatim}
7 7 1: 1 : 1
7 7 2: 1 : 1
7 7 3: 1 : 1
7 7 4: 1 : 1
7 7 5: 1 : 1
7 7 6: 1 : 1
7 7 7: 1 1 : 2
7 7 8: 1 2 : 3
7 7 9: 1 3 : 4
7 7 10: 1 4 : 5
7 7 11: 1 5 : 6
7 7 12: 1 6 : 7
7 7 13: 1 7 : 8
7 7 14: 1 8 1 : 10
7 7 15: 1 9 3 : 13
7 7 16: 1 10 6 : 17
7 7 17: 1 11 10 : 22
7 7 18: 1 12 15 : 28
7 7 19: 1 13 21 : 35
\end{Verbatim}
\normalsize
This is the number
of ways of placing $7\times 7$ squares into other squares:

\small
\begin{Verbatim}
 7 1 1:  1 : 1
 7 2 2:  1 : 1
 7 3 3:  1 : 1
 7 4 4:  1 : 1
 7 5 5:  1 : 1
 7 6 6:  1 : 1
 7 7 7:  1 1 : 2
 7 8 8:  1 4 : 5
 7 9 9:  1 9 : 10
 7 10 10:  1 16 0 : 17
 7 11 11:  1 25 0 : 26
 7 12 12:  1 36 0 : 37
 7 13 13:  1 49 0 0 : 50
 7 14 14:  1 64 126 28 1 : 220
 7 15 15:  1 81 468 440 79 : 1069
 7 16 16:  1 100 1128 2576 1246 0 : 5051
 7 17 17:  1 121 2220 9600 9550 0 : 21492
 7 18 18:  1 144 3870 27500 48175 0 0 : 79690
 7 19 19:  1 169 6216 66248 184681 0 0 0 : 257315
 7 20 20:  1 196 9408 141120 582904 0 0 0 0 : 733629
 7 21 21:  1 225 13608 277424 1719277 970368 166278 10062 207 1 : 3157451
 7 22 22:  1 256 18990 510220 4673549 8821336 5628634 1377284 129952 3600 : 21163822
\end{Verbatim}
\normalsize

\subsection{$8\times 8$ and larger Squares}
This is the number
of ways of placing $8\times 8$ squares into other squares:

\small
\begin{Verbatim}
 8 1 1:  1 : 1
 8 2 2:  1 : 1
 8 3 3:  1 : 1
 8 4 4:  1 : 1
 8 5 5:  1 : 1
 8 6 6:  1 : 1
 8 7 7:  1 : 1
 8 8 8:  1 1 : 2
 8 9 9:  1 4 : 5
 8 10 10:  1 9 : 10
 8 11 11:  1 16 : 17
 8 12 12:  1 25 0 : 26
 8 13 13:  1 36 0 : 37
 8 14 14:  1 49 0 0 : 50
 8 15 15:  1 64 0 0 : 65
 8 16 16:  1 81 160 32 1 : 275
 8 17 17:  1 100 582 500 79 : 1262
 8 18 18:  1 121 1380 2912 1246 0 : 5660
 8 19 19:  1 144 2680 10800 9550 0 : 23175
 8 20 20:  1 169 4620 30800 48175 0 0 : 83765
 8 21 21:  1 196 7350 73892 184681 0 0 : 266120
 8 22 22:  1 225 11032 156800 582904 0 0 0 : 750962
 8 23 23:  1 256 15840 303552 1591416 0 0 0 0 : 1911065
 8 24 24:  1 289 21960 552048 4113256 2087456 309736 15548 261 1 : 7100556
 8 25 25:  1 324 29590 952868 10029331 17541392 9842286 2039328 160080 3600 : 40598800
\end{Verbatim}
\normalsize

This is the number
of ways of placing $9\times 9$ squares into other squares:

\small
\begin{Verbatim}
 9 1 1:  1 : 1
 9 2 2:  1 : 1
 9 3 3:  1 : 1
 9 4 4:  1 : 1
 9 5 5:  1 : 1
 9 6 6:  1 : 1
 9 7 7:  1 : 1
 9 8 8:  1 : 1
 9 9 9:  1 1 : 2
 9 10 10:  1 4 : 5
 9 11 11:  1 9 : 10
 9 12 12:  1 16 : 17
 9 13 13:  1 25 0 : 26
 9 14 14:  1 36 0 : 37
 9 15 15:  1 49 0 : 50
 9 16 16:  1 64 0 0 : 65
 9 17 17:  1 81 0 0 : 82
 9 18 18:  1 100 198 36 1 : 336
 9 19 19:  1 121 708 560 79 : 1469
 9 20 20:  1 144 1656 3248 1246 : 6295
 9 21 21:  1 169 3180 12000 9550 0 : 24900
 9 22 22:  1 196 5430 34100 48175 0 : 87902
 9 23 23:  1 225 8568 81536 184681 0 0 : 275011
 9 24 24:  1 256 12768 172480 582904 0 0 0 : 768409
 9 25 25:  1 289 18216 332928 1591416 0 0 0 : 1942850
 9 26 26:  1 324 25110 598500 3882825 0 0 0 0 : 4506760
 9 27 27:  1 361 33660 1023300 9050334 4161822 542652 22964 321 1 : 14835416
\end{Verbatim}
\normalsize

This is the number
of ways of placing $10\times 10$ squares into other squares:

\small
\begin{Verbatim}
 10 1 1:  1 : 1
 10 2 2:  1 : 1
 10 3 3:  1 : 1
 10 4 4:  1 : 1
 10 5 5:  1 : 1
 10 6 6:  1 : 1
 10 7 7:  1 : 1
 10 8 8:  1 : 1
 10 9 9:  1 : 1
 10 10 10:  1 1 : 2
 10 11 11:  1 4 : 5
 10 12 12:  1 9 : 10
 10 13 13:  1 16 : 17
 10 14 14:  1 25 : 26
 10 15 15:  1 36 0 : 37
 10 16 16:  1 49 0 : 50
 10 17 17:  1 64 0 : 65
 10 18 18:  1 81 0 0 : 82
 10 19 19:  1 100 0 0 : 101
 10 20 20:  1 121 240 40 1 : 403
 10 21 21:  1 144 846 620 79 : 1690
 10 22 22:  1 169 1956 3584 1246 : 6956
 10 23 23:  1 196 3720 13200 9550 0 : 26667
 10 24 24:  1 225 6300 37400 48175 0 : 92101
 10 25 25:  1 256 9870 89180 184681 0 0 : 283988
 10 26 26:  1 289 14616 188160 582904 0 0 : 785970
 10 27 27:  1 324 20736 362304 1591416 0 0 0 : 1974781
 10 28 28:  1 361 28440 649800 3882825 0 0 0 : 4561427
 10 29 29:  1 400 37950 1101100 8660575 0 0 0 0 : 9800026
 10 30 30:  1 441 49500 1790580 18578835 7797510 904530 32706 387 1 : 29154491
\end{Verbatim}
\normalsize

The data base constructed above allows some extrapolations
while $s>1$:
\begin{conj}
\begin{equation}
T_{2s\times (2s+1)}(s,2)=1+10s+4s^2.
\end{equation}
\begin{equation}
T_{2s\times (2s+1)}(s,3)=2+16s.
\end{equation}
\begin{equation}
T_{2s\times (2s+1)}(s,4)=9.
\end{equation}
\end{conj}
While $s>2$:
\begin{conj}
\begin{equation}
T_{2s\times (2s+2)}(s,2)=3+18s+7s^2.
\end{equation}
\begin{equation}
T_{2s\times (2s+2)}(s,3)=8+40s.
\end{equation}
\begin{equation}
T_{2s\times (2s+2)}(s,4)=36.
\end{equation}
\end{conj}

\clearpage
\appendix
\section{Java Program for a Maple Generator}
The two Java programs \texttt{Hei.java} and \texttt{Tmat.java} 
in the \texttt{anc} subdirectory emit a Maple program that in principle generates
the inverse of $\mathbb{1}-M(z,t)$, where $\mathbb{1}$ is the
unit matrix and $M$ the transfer matrix.  Since we are only interested in
the top left element of the inverse---which is the bivariate generating function---,
effectively only a linear system of equations is set up.
The programs are compiled with
\begin{verbatim}
javac -cp . *.java
\end{verbatim}
and the output of the main program of \texttt{Tmat} can be directly
piped into Maple:
\verb+ java -cp . Tmat +\textit{s} \textit{n} \verb+|maple -q +

The two command line options $s$ and $n$
are the edge length of the squares and the width of
the rectangle.
The programs are licensed under the LGPLv3.

\bibliographystyle{amsplain}
\bibliography{all}

\end{document}